\documentclass[nofootinbib,showpacs,showkeys]{revtex4}%
\usepackage{etoolbox}
\usepackage[colorlinks=true]{hyperref}
\AtBeginDocument{
	\hypersetup{
		urlcolor=cyan,
		citecolor=magenta,
		linkcolor=blue,
		anchorcolor=green,
	}%
}
\usepackage{amssymb}
\usepackage{amsmath}
\usepackage{extpfeil}
\usepackage{mathtools}
\usepackage{amsfonts}
\usepackage{graphicx}
\usepackage{ulem} 
\usepackage{cancel}
\usepackage[usenames]{xcolor}
\usepackage{epstopdf}
\usepackage{extarrows}
\usepackage{multirow}

\usepackage{tikz}
\usetikzlibrary{matrix,arrows}

\def\mod#1{\,({\rm mod\ }#1)}
\def\p{\alpha}

\definecolor{amaranth}{rgb}{0.9, 0.17, 0.31}
\definecolor{purple(munsell)}{rgb}{0.62, 0.0, 0.77}
\definecolor{americanrose}{rgb}{1.0, 0.01, 0.24}
\definecolor{palatinateblue}{rgb}{0.15, 0.23, 0.89}
\definecolor{royalblue(web)}{rgb}{0.25, 0.41, 0.88}
\definecolor{hanpurple}{rgb}{0.32, 0.09, 0.98}
\definecolor{beaublue}{rgb}{0.74, 0.83, 0.9}
\definecolor{carminered}{rgb}{1.0, 0.0, 0.22}
\definecolor{brightpink}{rgb}{1.0, 0.0, 0.5}
\definecolor{vividviolet}{rgb}{0.62, 0.0, 1.0}
\hypersetup{ linktoc=all,
	colorlinks, linkcolor={palatinateblue},
	citecolor={brightpink}, urlcolor={amaranth}}

\def\sideremark#1{\ifvmode\leavevmode\fi\vadjust{\vbox to0pt{\vss
			\hbox to 0pt{\hskip\hsize\hskip1em
				\vbox{\hsize2cm\tiny\raggedright\pretolerance10000
					\noindent #1\hfill}\hss}\vbox to8pt{\vfil}\vss}}}%

\makeatletter
\newcommand{\qbar}{\text{\q@bar}}
\newcommand{\q@bar}{%
	\vphantom{$\m@th q$}%
	\ooalign{%
		$\m@th q$\cr
		\hidewidth\kern.15em\smash{\raisebox{-0.7ex}{$\m@th\mathchar'26$}}\hidewidth\cr}%
}
\makeatother

\begin{document}

\title[ ]{Cauchy distributions for the integrable standard map  }
\author{Anastasios Bountis}
\email{anastasios.bountis@nu.edu.kz}
\affiliation{Department of Mathematics,
	Nazarbayev University, Nur--Sultan, Kazakhstan}
\author{ J. J. P. Veerman}
\email{veerman@pdx.edu}
\affiliation{$^{2}$Department of Mathematics and Statistics§,
	Portland State University, Portland, OR, USA}
\author{Franco Vivaldi}
\email{f.vivaldi@maths.qmul.ac.uk}
\affiliation{$^{3}$School of Mathematical Sciences, Queen Mary,
	University of London,
	London E1 4NS, UK}


\begin{abstract}
We consider the integrable (zero perturbation) two--dimensional standard map, in light
of current developments on ergodic sums of irrational rotations, and recent
numerical evidence that it might possess non-trivial $q$-Gaussian statistics. 
Using both classical and recent results, we show that
the phase average of the sum of centered positions of an
orbit, for long times and after normalization, obeys the Cauchy distribution
(a $q$-Gaussian with $q=2$), while for almost all individual orbits such a sum
does not obey any distribution at all. We discuss the question of existence of distributions for KAM tori.
\end{abstract}

\eid{ }
\date{\today }
\startpage{1}
\endpage{1}
\maketitle

\section{Introduction}\label{Introduction}
In the last fifteen years, there has been a growing number of numerical investigations suggesting that chaotic orbits of conservative (primarily Hamiltonian) systems are characterized by $q$--Gaussian statistics, whose index varies from $q=1$ (Gaussian) in `wide chaotic seas' to $1<q<3$ in `thin chaotic layers' close to invariant tori (see \cite[ch. 8]{BouSko2012} for details).

More recently, numerical evidence has shown that, as the perturbation
parameter $\varepsilon>0$ decreases, chaotic orbit distributions
of the Chirikov-Taylor standard map (a paradigm of two--degree--of--freedom conservative flows)
\begin{equation}\label{eq:StandardMap}
(x,y)\mapsto(x+y,y+\varepsilon\sin(2\pi x)) \mod{1}
\end{equation}
change from Boltzmann-Gibbs (Gaussian) to Tsallis ($q$-Gaussian) statistics
\cite{TirnakliBorges}. At present, for $\varepsilon\not=0$, this
claim cannot be rigorously verified nor challenged.

The purpose of this Letter is to elucidate the ergodic properties of
the standard map for $\varepsilon=0$ by applying the theory of ergodic
sums of irrational rotations. This contribution fits within the
broader study of regular motions on two-dimensional invariant tori
and its connection with their (irrational) rotation number
(see \cite{LictLieb1992} and the collection of articles in \cite{MacKayMeiss1987}).

For $\varepsilon=0$ the momentum $y$ in (\ref{eq:StandardMap}) is constant,
the map is integrable, and the phase space foliates into invariant circles (tori),
parametrized by $y$. Let $y=\alpha$ represent one such circle, with $\alpha\not\in\mathbb{Q}$, and
consider the distributional properties (if any) of the ergodic sums
\begin{equation}\label{eq:S}
S_t(x,\p)=\sum_{k=0}^{t-1} h(f^k_\p(x)),\qquad
h(x)=\{x\}-\frac{1}{2},\qquad
f_\p(x)=\{x+\p\}
\end{equation}
where $\{\cdot\}$ denotes the fractional part, within the interval $[0,1)$ with end-points identified.

We will show that the distributional properties of (\ref{eq:S}) for
$\varepsilon=0$ can be rigorously determined by applying both classical
and recent results of ergodic theory that deal with
number-theoretic questions.
These are theorems that focus on the convergence in distribution of the sums
$S_t(x,\p)$, \textit{appropriately centralized and normalized}.
As $t\to\infty$, our main findings are (details will be given below):
{\it
	\begin{enumerate}
		\item [1.] The phase average of $S_t$ converges to a Cauchy
		distribution (a $q$-Gaussian with $q=2$).
		\item [2.] For almost all initial points $(x,\p)$,
		$S_t$ does not converge to any distribution.
		\item [3.] If $\p$ is a badly approximable irrational\footnote{with bounded
			continued fraction coefficients}, then $S_t$ converges in distribution
		to a Gaussian (central limit theorem).
	\end{enumerate}
}
The first statement follows from a classical theorem by Kesten (1960); 
the second is due to Dolgopyat \& Sarig (2018); 
the third to Bromberg and Ulcigrai (2018). 

Based on the above, we shall proceed to examine recent numerical evidence by Tirnakli \& Tsallis \cite{TirnakliTsallis} on the sum (\ref{eq:S}) for the map (\ref{eq:StandardMap}) with $\epsilon=0$. More specifically, they sought to identify its statistical distribution within the family of $q$-Gaussians with zero mean
\begin{equation}\label{eq:q-Gaussian}
P_q(s)=\frac{\sqrt{\beta}}{C_q}\exp_q(-\beta s^2)\quad\mbox{where}\quad
\exp_q(x)=[1+(1-q)x]^{\frac{1}{1-q}}
\end{equation}
where $\beta$ is a parameter and $C_q$ is the normalization coefficient \cite{UmarovEtAl}.
The Cauchy distribution of Theorem 1 that we obtain for this case corresponds to $q=2$, with $C_2=\pi$.
In spite of approximately $8\times 10^{14}$ iterations of the map ($2\times10^8$ initial points, each iterated $2^{22}$ times), the numerical value proposed in \cite{TirnakliTsallis} $q\approx1.935$ agrees with our theoretical prediction, within only one significant digit (see Appendix A).

\section{Background theorems}\label{Theorems}
\medskip

\begin{figure}[t]
	\hspace*{-20pt}
	\begin{minipage}{6cm}
		\centering
		\includegraphics[scale=0.42]{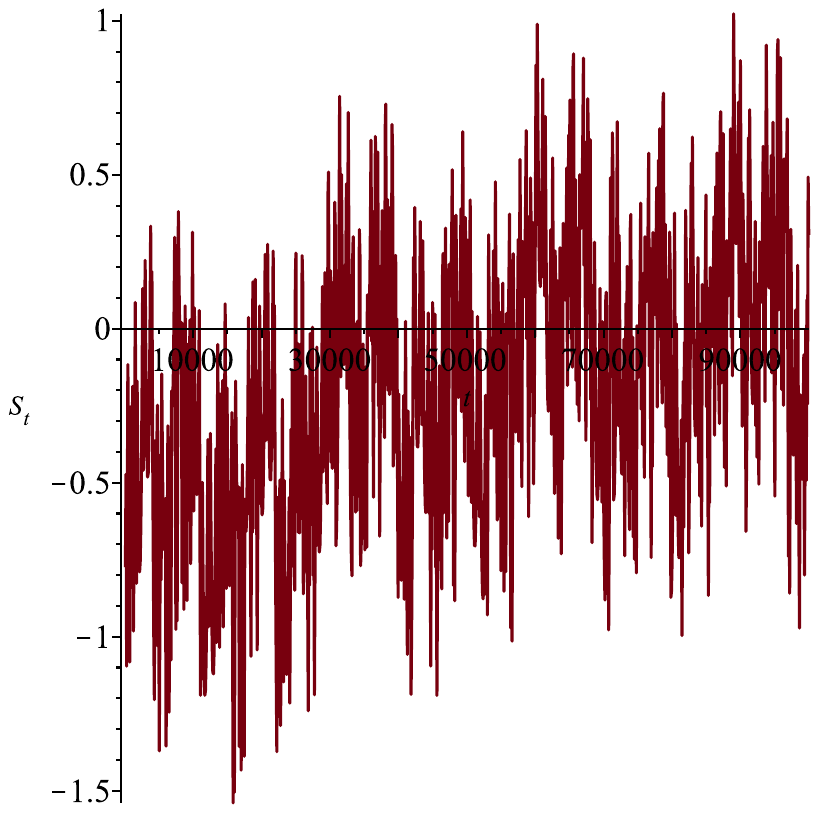}\\
	\end{minipage}
	\hspace*{40pt}
	\begin{minipage}{6cm}
		\centering
		\includegraphics[scale=0.42]{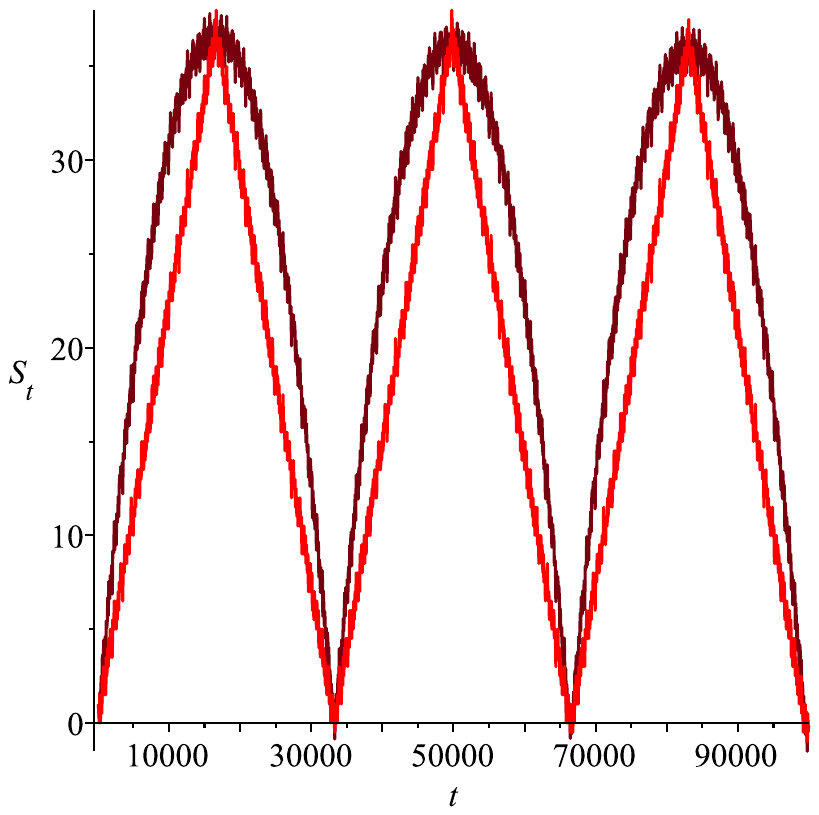}\\
	\end{minipage}
	\hspace*{50pt}
	\hfil
	\vspace*{-140pt}
	\caption{\rm\small
		Left: The sequence $S_t$, with initial condition $x=0$ and $t\leqslant 10^5$.
		Left: $\p=(\sqrt{5}-1)/2$, with continued fractions coefficients
		$[1,1,1,\ldots]$.
		Right: $\p=\pi-3$ with continued fractions $[7,15,1,293,1,\ldots]$.
		The isolated large coefficient causes large fluctuations with an approximate
		periodicity of 33102, the denominator of the corresponding convergent.
		The red curve corresponds to the observable $-\mathrm{sign}(h(x))/2$, see (\ref{eq:h2}).
	}
	\label{fig:S}
\end{figure}

Let us begin with some general observations:

\noindent 1.\/
The distributional properties of $S_t$ depend sensitively on the type of averaging performed. Indeed, in (\ref{eq:S}) one may randomise any combination of $x,\p$,
and $t$, each sampled from a uniform distribution. In this regard, there are several
types of limit theorems, stated below.

\noindent 2.\/
The behaviour of $S_t$ for individual orbits is dictated by the arithmetical
properties of $\p$.  Bounded continued fraction coefficients improve convergence
(Figure \ref{fig:S}, left), while the presence of large coefficients in the continued fraction expansion
cause large fluctuations (Figure \ref{fig:S}, right).
Recall that almost all real numbers have unbounded coefficients \cite[Theorem 29]{Khintchin}.

\noindent 3.\/
The properties of $S_t$ in (\ref{eq:S}) depend on the choice of the observable $h$.
Thus, if $h$ is sufficiently smooth, then $S_t$ is bounded for almost all pairs $(x,p)$
\cite{Herman} (see also \cite[Appendix A]{DolgopyatSarig:20}). Bounded ergodic sums may
occur also with non-smooth observables \cite{Kesten:66}, but not for our choice of $h$.
\medskip

Let us now recall various limit theorems in \cite{DolgopyatSarig:17,DolgopyatSarig:20}, adopting the corresponding terminology.

For distributions along individual orbits, we fix the initial
conditions and randomise $t$.
The ergodic sum $S_t$ of a map $f$ and observable $h$ satisfies a
\textit{temporal distributional limit theorem} if there is
a centralizing sequence $U_t$, a normalizing sequence
$V_t\to\infty$, and a (non-constant) random variable $Y$ such that
for all real $y$
\begin{equation}\label{eq:TDLT}
\lim_{T\to\infty}\frac{1}{T}\#\Bigl\{t\in\{0,\ldots,T-1\}\,:\,
\frac{S_t(x,\p)-U_T}{V_T}\leqslant y\Bigr\}=F_Y(y),
\end{equation}
where \# denotes the cardinality and $F_Y$ is the cumulative distribution function of $Y$.
Equivalently, the random variable $(S_t-U_T)/V_T$, where $t$ is
uniformly distributed among the first $T$ iterates,
converges in distribution to $Y$, as $T$ goes to infinity.
In general, the quantities $U_T,V_T$, and $Y$ will depend on the
initial point, as well as on the function $h$. 

In a \textit{spatial distributional limit theorem} one randomizes $x$ instead of $t$,
and modifies (\ref{eq:TDLT}) as follows:
\begin{equation}\label{eq:SDLT}
\lim_{T\to\infty}\frac{1}{T}\,\mu\Bigl\{x\in [0,1)\,:\,\frac{S_T(x,\p)-U_T}{V_T}\leqslant y\Bigr\}=
F_Y(y),
\end{equation}
where $\mu$ is the Lebesgue measure. (In a more general setting,
the unit interval is replaced by the phase space of $f$,
with invariant measure $\mu$.)
Unlike in (\ref{eq:TDLT}), here only $S_T$ is considered, the earlier
values of the sum being ignored.

In the above limit theorems $\p$ is kept fixed, being regarded as a parameter.
If we randomize $\p$ then the limit theorems (spatial or temporal) are said to
be \textit{annealed}. For the integrable standard map, a limit theorem
resulting from a phase average is of the annealed spatial type.


Let us recall the first limit theorem of this kind, regarding rotations, due to Kesten \cite{Kesten:60}:
\medskip

\noindent\textbf{Theorem 1} \cite{Kesten:60}.
\textit{If $(x,\p)$ is uniformly distributed on $\mathbb{T}^2$,
	then the distribution of $S_T(x,\p)/\ln T$ converges as
	$T\to\infty$ to a Cauchy distribution: for some $\rho$ and all $y$
	we have
	\begin{equation}\label{eq:Cauchy}
	\lim_{T\to\infty}\frac{1}{T}\mu\Bigl\{(x,\p)\in\mathbb{T}^2\,:\,\frac{S_T(x,\p)}{\ln T}\leqslant y\Bigr\}
	=\frac{\rho}{\pi}\int_{-\infty}^{y}\frac{\mathrm{d}x}{1+\rho^2x^2},
	\end{equation}
	where $\mu$ is the Lebesgue measure.
}
\medskip

Thus, phase averaging and logarithmic scaling yield the Cauchy distribution.
Kesten, in fact, gives a formula for the value of $\rho$, which we compute
in Appendix B below to find $\rho=4\pi$.
Comparing the right-hand-side of (\ref{eq:Cauchy})
with (\ref{eq:q-Gaussian}) for $q=2$, we find $\beta=\rho^{2}$, hence
$P(0)=\rho/\pi=4$.

Recently, Dolgopyat \& Sarig \cite[Theorem 2.1]{DolgopyatSarig:20}
established a temporal version of Kesten's theorem, randomising
$t$ and $\p$ instead of $x$ and $\p$.
The distribution is again Cauchy, with the same logarithmic
scaling but a different constant: $\rho=3\pi\sqrt{3}$.

\medskip

We now turn to the existence of temporal distributions for individual orbits.
The initial condition $(x,\p)$ is fixed, with $\p\not\in\mathbb{Q}$,
and we attempt to extract a distribution from the terms of the sequence (\ref{eq:S}).
The following result (Dolgopyat \& Sarig, 2018) shows that without $\p$-averaging, almost surely,
no temporal distribution exists.
\medskip

\noindent\textbf{Theorem 2} \cite[Theorem 1.2]{DolgopyatSarig:18}.
\textit{Let $h$ be a piecewise smooth function of zero mean.
	Then there is a set of full measure $\Lambda\subset\mathbb{T}^2$
	such that, if $(x,\p)\in\Lambda$ then the ergodic sum of $h\circ f_\p$
	does not satisfy a temporal distributional limit theorem
	on the orbit of $x$.}
\medskip

The observable $h$ in (\ref{eq:S}) satisfies the assumptions of this theorem.
Thus, for a generic initial condition $(x,\p)$, the sequence (\ref{eq:S})
for the integrable standard map does not admit any distributional limit,
which is our second statement in Section \ref{Introduction}.
This is because there are different scaling limits on different
subsequences, due to the presence of large continued fractions coefficients
The latter appear at random in the continued fraction expansion of
most numbers, due to the ergodic properties of Gauss' map
(see \cite[section 3.2]{EinsiedlerWard} and
\cite[Section 1.3]{DolgopyatSarig:18}).
Thus, distributional data extracted from a typical orbit are
intrinsically unstable, and this problem cannot be fixed by choosing
appropriate centralizing and normalizing sequences in (\ref{eq:TDLT}).

The capricious nature of the ergodic sums (\ref{eq:S}) is best illustrated 
by fixing $\alpha$ and randomizing the initial condition. For any fixed time $t$, 
the appropriately normalized $S_t$ of \eqref{eq:S} gives a well-defined 
density $\nu(z,t,\alpha)$, symmetric around the origin (see \cite{VF}).

\begin{figure}[pbth]
	\begin{center}
		\includegraphics[height=4.3cm]{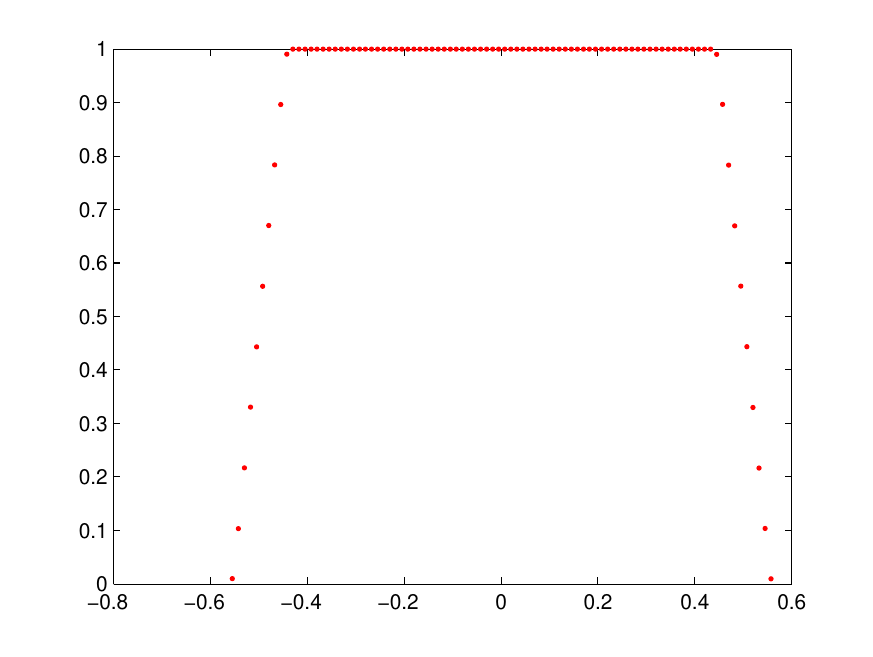}
		\includegraphics[height=4.3cm]{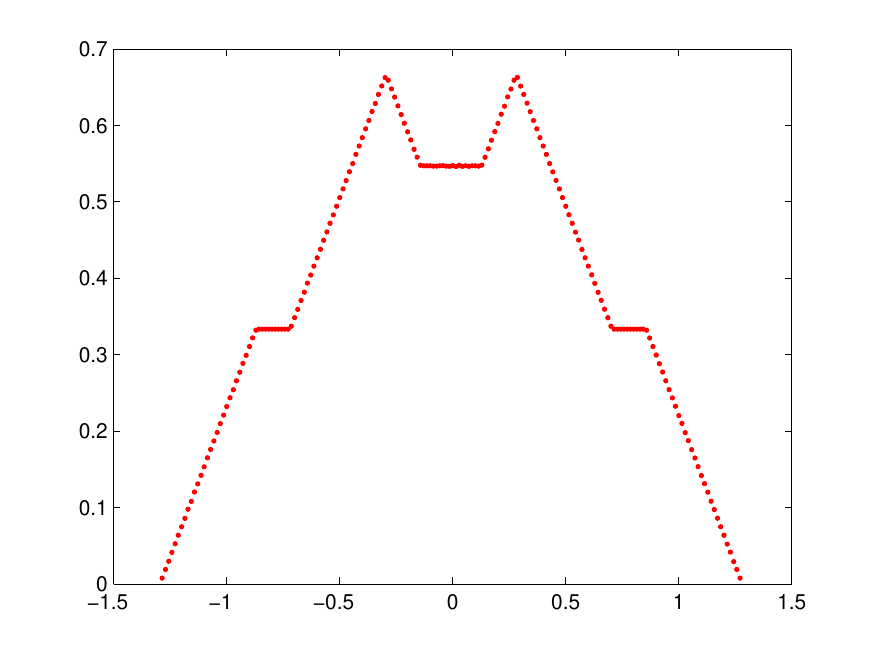}
		\includegraphics[height=4.3cm]{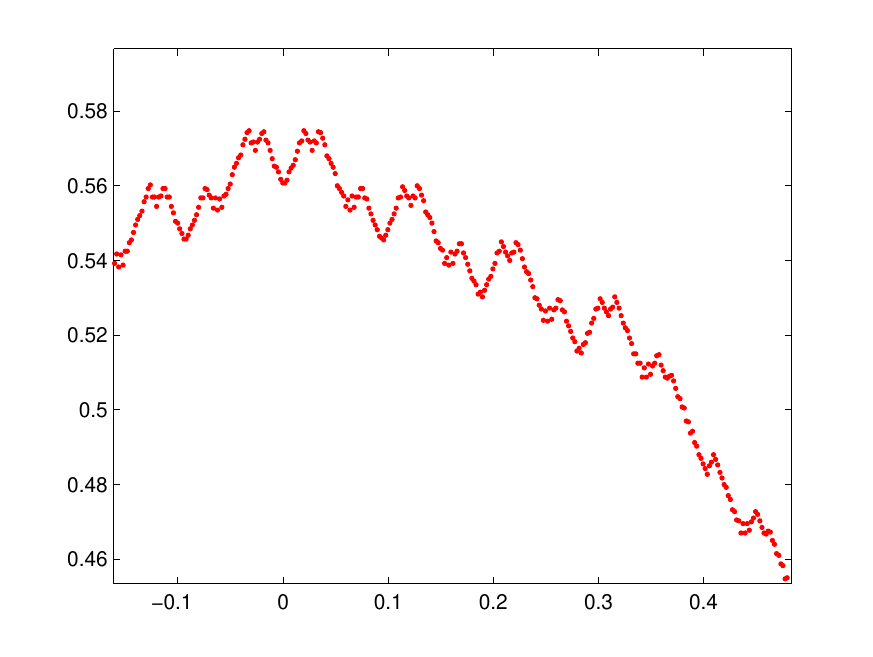}
		\caption{\rm\small Numerical approximations of the density $\nu(z,t,\alpha)$
			for $\alpha=\mathrm{e}-2$, and $t=1001$ (left), $t=213$ (centre), and
			a detail for $t=334$ (right).
			Only the first value of $t$ is a denominator of a convergent of $\mathrm{e}-2$.}
		\label{fig:emintwo}
	\end{center}
\end{figure}

In Figure \ref{fig:emintwo} we plot the density $\nu(z,t,\alpha)$ as a function 
of $z$, for $\alpha=\textrm{e}-2=[1,2,1,1,4,1,1,6,\ldots]$, with unbounded 
(although still regular) continued fractions coefficients. As $t$ changes, we 
observe significant variations, ranging from near-uniformity, if $t$ is the 
denominator of a convergent, to exotic shapes for other values of $t$.
For a typical $\alpha$ these variations cannot be tamed by an averaging process.


\vskip -0.0cm

To obtain temporal limit theorems for individual orbits, one must restrict $\p$ to a zero measure set over
which the fluctuations of $S_t$ can be controlled.
In this direction, a \textit{temporal central limit theorem} has been proved very recently by Bromberg and
Ulcigrai (see below) who considered the sum $S_t$ corresponding to the much--studied
observable \cite{Kesten:60,Kesten:62,Kesten:66,DrmotaTichy}
\begin{equation}\label{eq:h2}
h_2(x)=\chi_I(x)-\gamma,
\end{equation}
where $I=[0,\gamma)$ is an interval of length $\gamma$ and $\chi_I$ is the
characteristic function of $I$.
(For $\gamma=1/2$ we have $h_2(x)=-\mathrm{sign}(h(x))/2$, see Figure \ref{fig:S}, right.)
The zero--measure set chosen here consists of the \textit{badly approximable}
values of $\p$,  i.e., those with bounded continued fraction coefficients.
This set is `large' in that it has Hausdorff dimension 1 \cite{Jarnik}.
\medskip

\noindent\textbf{Theorem 3} \cite[Theorem 1.1]{BrombergUlcigrai}.
{\it Let $\p$ be a badly approximable irrational number. For every $\gamma$
	badly approximable with respect to $\p$, every $x$, and every real $y$ we have
	$$
	\frac{1}{T}\#\Bigl\{0\leqslant t< T\,:\,\frac{S_t(\alpha,\gamma,x)-U_t}{V_t}\leqslant y\Bigr\}
	\to\frac{1}{\sqrt{2\pi}}\int_{-\infty}^ye^{-x^2/2}\mathrm{d}x
	$$
	for some sequences $U_t(\alpha,x,\gamma)$ and $V_t(\alpha,\gamma)\to\infty$.
}
\medskip

The restriction on $\gamma$ is a (mild) diophantine condition
(see \cite{BrombergUlcigrai} for details).
This theorem is our third statement in Section \ref{Introduction}. It generalizes previous results by
Beck \cite{Beck1,Beck2}, who dealt with the special case of quadratic
irrational $\alpha$ (irrational roots
of a quadratic polynomial with integer coefficients), rational $\gamma$,
and $x=0$. Under these constraints, he was able to obtain a good description of
centralizing and normalizing sequences: $U_t=U\ln t$ and $V_t=V\sqrt{\ln t}$,
for some constants $U$ and $V$.
\medskip

\section{Concluding remarks}\label{Conclusions}

In this paper we have considered the problem of the statistical properties
of regular motions in the zero perturbation limit of the two--dimensional standard map.
We remark at the outset, regarding the sensitivity of ergodic sums on the choice of the
observable, that great care must be taken when assessing the relevance
of the theory reviewed above to more general settings
(cf.~ \cite[section 2]{DolgopyatSarig:18}).

Our main finding is that orbits on tori of the integrable two--dimensional standard map are
described by the Cauchy distribution. We conjecture that this may hold more generally for integrable conservative systems. A second conclusion that follows from our findings is that the almost certain non-existence of distributions for individual orbits that holds for a large class of observables
\cite[section 1]{DolgopyatSarig:20} is expected to apply to
KAM tori as well.

Regarding the existence and type of distributions resulting from averaging
over KAM tori, we note that the latter are parametrized by a positive measure
Cantor set of rotation numbers, selected from the unit interval via a
perturbation-dependent diophantine condition \cite[p.~344]{ArrowsmithPlace}.
This condition removes from the averaging set the rotation numbers that are
too closely approximable by rationals, thereby avoiding large
fluctuations of the ergodic sum and improving convergence.

Thus, it seems plausible that, for appropriate observables, an ergodic sum
such as (\ref{eq:S}), averaged over KAM tori, will converge in distribution.
If so, then what would be the resulting distribution as a function of the
perturbation parameter? This question appears worthy of further study.

\begin{acknowledgments}
 We are grateful to U. Tirnakli and C. Tsallis for useful discussions.
\end{acknowledgments}

\bigskip

\appendix
\section{The results of Tirnakli and Tsallis}

Let us examine the results of the numerical experiments presented in
\cite{TirnakliTsallis}, in view of Theorem 1 of Section \ref{Theorems}.
The authors choose a large set of $2\times 10^8$ initial points
uniformly distributed on the 2-torus and compute
$T=2^{22}\approx 4\times 10^6$ iterates of (\ref{eq:StandardMap})
for each initial point, letting $\langle x\rangle$ be the average
of the resulting set of $8\times 10^{14}$ data points.
For each orbit, the final value $S_T$ of the ergodic sum is stored,
and centralized using $\langle x\rangle$.
The $2\times 10^8$ values of $S_T$ thus obtained are merged together,
their distribution is computed and scaled by the numerical value of
$P(0)$, while the latter, together with the normalisation coefficient $C_q$
\cite{UmarovEtAl} provide a numerical value for $\beta$.
Fitting the normalized distribution to their data yields the value
$q\approx1.935$.

The observable used in \cite{TirnakliTsallis} differs from that of
Theorem 1, since $1/2$ is replaced by $\langle x\rangle$, although
such a difference can be absorbed by a centralized sequence $U_T$.
If we apply Theorem 1 to these data, the numerical result for
$q$ agrees within only one significant figure. The numerical value of P(0) 
(after employing $ln(T)$ scaling) is $P(0)\approx 1.5$ \cite{private}, 
which is also problematic when compared with the exact value we have obtained, $P(0)=4$.
In order to investigate these discrepancies further, at this stage, 
we would need more detailed information about the computations presented
in \cite{TirnakliTsallis}.
\medskip

\section{Computation of $\rho$ in Theorem 1}

For $x\in [0,1)$, let $q_n$ be the denominator of the $n$th convergent
of the continued fractions of $x$ \cite[Chapter X]{HardyWright}.
In \cite[p. 66]{Khintchin}
it is proved that
there is a unique $\tau > 0$ such that for almost all $x\in(0,1)$ the
following limit exists
\begin{equation*}
\tau=\lim_{n\rightarrow \infty} \dfrac{n}{\ln q_n(x)} =\dfrac{12 \ln 2}{\pi^2}.
\end{equation*}
The parameter $\rho$ is given by the following formula
\cite{Kesten:60}
\begin{equation}\label{eq:rho}
\rho=\frac{2\pi\ln 2}{\tau I}
\end{equation}
where
\begin{equation}\label{eq:I}
I=\int_0^1 \int_0^1\,\left|\,\sum_{k=1}^\infty\,k^{-2}\,\sin 2\pi k x
\sin 2\pi k y\right|\,\mathrm{d}y\mathrm{d}x.
\end{equation}

The integrand has an 8-fold symmetry, being invariant under reflection
with respect to the main diagonal and each of the lines $x=1/2,y=1/2$
(because of the absolute value).
So it suffices to restrict the integration to the triangle with
vertices $(0,0),(1,0),(1/2,1/2)$. The sum is absolutely convergent,
and over that domain can be shown (via Fourier analysis) to be
equal to $\pi^2 y(1-2x)$, which is positive.
So, by Fubini's theorem, we have
$$
I=8\int_0^{1/2}\mathrm{d}x\int_0^x\pi^2 y(1-2x)\mathrm{d}y
=\frac{\pi^2}{24}
$$
and thererefore $\rho=4\pi$ as desired.

\medskip


\end{document}